\documentclass[12pt]{article}
\usepackage{amssymb,amsmath,amsthm,latexsym}

\newtheorem{thm}{Theorem}[section]

\newtheorem{lem}[thm]{Lemma}
\newtheorem{cor}[thm]{Corollary}
\newtheorem{exam}{Example}
\newtheorem{defi}[thm]{Definition}
\newtheorem{rem}[thm]{Remark}
\newcommand{\pf}{{\bf Proof. \ }}

\begin{document}

\title{ On the arithmetic of  the endomorphism ring End($\mathbb{Z}_{p}\times\mathbb{Z}_{p^{m}}$)}
 \author{
 Xiusheng Liu\\
 School of Mathematics and Physics, \\
 Hubei Polytechnic University  \\
 Huangshi, Hubei 435003, China, \\
{Email: \tt lxs6682@163.com} \\
Hualu Liu\\
 School of Mathematics and Physics, \\
 Hubei Polytechnic University  \\
 Huangshi, Hubei 435003, China, \\
{Email: \tt hwlulu@163.com} \\}
\maketitle



\begin{abstract}
For a prime $p$, let $E_{p,p^m}=\{\begin{pmatrix}a&b\\p^{m-1}c&d\end{pmatrix}|a,b,c\in\mathbb{Z}_{p},~\mathrm{and}~d\in \mathbb{Z}_{p^{m}}\}$. We first establish a ring isomorphism from $\mathrm{End}(\mathbb{Z}_p\times\mathbb{Z}_p^m)$ onto $E_{p,p^m}$. We then provide the way to compute $-d$ and $d^{-1}$ using arithmetic in  $\mathbb{Z}_{p}$ and $\mathbb{Z}_{p^{m}}$, and characterize invertible elements in $E_{p,p^m}$. Moreover, we introduce the minimal polynomial for each element in $E_{p,p^m}$ and given its applications.

\end{abstract}

\bf Keywords\rm : Endomorphism ring$\cdot$Invertible element $\cdot$Minimal polynomials

\bf 2000 MR Subject Classification\rm:  16S50$\cdot$16W20$\cdot$13B25

\section{Introduction}
 Throughout this paper, $p$ denotes a prime number.  Bergman[2] established the $\mathrm{End}(\mathbb{Z}_{p} \times \mathbb{Z}_{p^{2}})$   is a semilocal ring with $p^3$ elements that cannot be embedded in matrices over any commutative ring. Climent et al. [1] identitied the elements of $\mathrm{End}(\mathbb{Z}_{p} \times \mathbb{Z}_{p^{2}})$ with elements in a new set, denoted by $E_p$, of  matrices of size $2\times2$, whose elements in the first row belong to $\mathbb{Z}_p$ and the elements in the second row belong to $\mathbb{Z}_{p^{2}}$; also, using the arithmetic in $\mathbb{Z}_{p}$ and $\mathbb{Z}_{p^{2}}$, then characterized the invertible elements of $E_p$.  In this paper, we propose an approach following [1] to investigate the arithmetic of the generalized of endomorphism rings  $\mathrm{End}(\mathbb{Z}_{p} \times \mathbb{Z}_{p^{m}})$.

 Recall that $\mathbb{Z}_n=\{0,1,2,\ldots,n-1\}$ is a commutative unitary ring with the addition and multiplication mod $n$, i.e.,  for any $a,b\in\mathbb{Z}_n$
 $$(a+b)~\mathrm{in} ~\mathbb{Z}_n=((a+b)~{in} ~\mathbb{Z})~\mathrm{mod}~n)$$
 and
 $$(ab)~\mathrm{in} ~\mathbb{Z}_n=((ab)~{in} ~\mathbb{Z})~\mathrm{mod}~n).$$

 From now on,  let we assume that $m$ is a positive integer and discuss the ring $\mathbb{Z}_p$ and $\mathbb{Z}_{p^m}$.  Clearly, we can also assume that $\mathbb{Z}_p\subset\mathbb{Z}_{p^m}$, even though  $\mathbb{Z}_p$ is not a subring of $\mathbb{Z}_{p^m}$. Then, it follows that notation is utmost important to prevent errors like the following. Suppose that $p=7$, then
 $$\mathbb{Z}_7=\{0,1,2,4,5,6\},~\mathrm{and}~\mathbb{Z}_{7^{2}}=\{0,1,2,\ldots,48\}.$$
 Note that $3,5\in \mathbb{Z}_7$ and $3+5=1\in \mathbb{Z}_7$;  but $3,5\in \mathbb{Z}_{7^2}$ equally. However when $3,5\in \mathbb{Z}_{7^2}$, $3+5=8\in \mathbb{Z}_{7^2}$. Obviously, $1\neq8$ in   $\mathbb{Z}_{7^2}$. Such error can be easily avoidable if we write, when necessary, $x ~\mathrm{mod}~p$ and $x~\mathrm{mod}~p^m$ to refer the element $x $ when $x\in \mathbb{Z}_p$ and $x\in \mathbb{Z}_{p^m}$, respectively. In this light, the above example could be written as
 $$(3~\mathrm{mod}~7)+(5~\mathrm{mod}~7)=1~\mathrm{mod}~7,$$
 whereas
 $$(3~\mathrm{mod}~7^2)+(5~\mathrm{mod}~7^2)=8~\mathrm{mod}~7^2.$$

 This paper is structured as follows. In section 2, we construct a ring isomorphism from $\mathrm{End}(\mathbb{Z}_{p} \times \mathbb{Z}_{p^{m}})$ onto $E_{p,p^m}$. In section 3, we first provide the way to compute $-d$ and $d^{-1}$ using arithmetic in  $\mathbb{Z}_p$ and $\mathbb{Z}_{p^{m}})$.  Then we characterize invertible elements in $E_{p,p^m}$.  Finally, in section 4, we introduce the minimal polynomial for each element in $E_{p,p^m}$ and consider its applications.
\section{A characterization of the ring $\mathrm{End}(\mathbb{Z}_{p} \times \mathbb{Z}_{p^{m}})$}
In this section, we denote by $1_{p}$ and $1_{p^m}$ the identity of $\mathbb{Z}_p$ and $\mathbb{Z}_p^m$, respectively. In order  to establish a ring isomorphism from $\mathrm{End}(\mathbb{Z}_p\times\mathbb{Z}_p^m)$ onto $E_{p,p^m}$, let us recall a well-known result from ring theory(See [3]).

The ring $\mathbb{Z}_{p^{m}}$ is a Galois ring of characteristic $p^{m}$. Each element $d\in\mathbb{Z}_{p^{m}}$ can be written uniquely as a finite sun
$$d=p^{m-1}u_{m-1}+\cdots+pu_{1}+u_{0},$$
where $u_i\in\mathbb{Z}_{p}$, for $i=0,1,\ldots, m-1$. The invertible elements in $\mathbb{Z}_{p^{m}}$ are precisely the $d$ for which $u_0\neq0$,  and $pd=0$ if and only if $d=p^{m-1}a$, i.e., $d$ is not invertible in $\mathbb{Z}_{p^m}$.

Since $\mathbb{Z}_{p} \times \mathbb{Z}_{p^{m}}$ is a  $\mathbb{Z}_{p^{m}}$-module, an element $(u,v)\in\mathbb{Z}_{p} \times \mathbb{Z}_{p^{m}}$ may be written uniquely as
$$r(1_p,0)+s(0,1_{p^m}),$$
where $r\in\mathbb{Z}_{p}$, and  $s\in\mathbb{Z}_{p^{m}}.$

Let $\alpha\in \mathrm{End} (\mathbb{Z}_{p} \times \mathbb{Z}_{p^{m}})$.  Since $\alpha$ is $\mathbb{Z}_{p^{m}}$-linear, we have
$$\alpha(r(1_p,0)+s(0,1_{p^m}))=r\alpha(1_p,0)+s\alpha(0,1_{p^m}),$$
for any $r\in\mathbb{Z}_{p}$, and  $s\in\mathbb{Z}_{p^{m}}$. Therefore, $\alpha$ is uniquely determined by two elements $\alpha(1_p,0)$ and $\alpha(0,1_{p^m})$ of $\mathbb{Z}_{p} \times \mathbb{Z}_{p^{m}}$.

Now, suppose that $\alpha(1_p,0)=(a,e)$ and  $\alpha(0,1_{p^m})=(b,d)$ with $a,b\in\mathbb{Z}_{p}$ and $e,d\in\mathbb{Z}_{p^{m}}$. Since $pa=p1_p=0$, we have $(0,pe)=(pa,pe)=p(\alpha(1_p,0))=\alpha(p1_p,p0)=(0,0)$, which implies $pe=0$ in $\mathbb{Z}_{p^{m}}$. Hence there exists a uniquely $c\in\mathbb{Z}_{p}$ such that $e=p^{m-1}c$. This means that
$$ \mathrm{End} (\mathbb{Z}_{p} \times \mathbb{Z}_{p^{m}})=\{\alpha|\alpha(1_p,0)=(a,p^{m-1}c),\alpha(0,1_{p^m})=(b,d)\},$$
where $a,b,c\in\mathbb{Z}_{p}$, and $d\in\mathbb{Z}_{p^{m}}$.

Obviously, $\mathrm{End} (\mathbb{Z}_{p} \times \mathbb{Z}_{p^{m}})$ is a nonccommutative unitary ring with the usual addition and composition of endomorphisms.

The proof in following Theorem is similar to that for [2] Theorem 3.
\begin{thm} If $p$ is a prime number, then the ring of endomorphism  $\mathrm{End} (\mathbb{Z}_{p} \times \mathbb{Z}_{p^{m}})$ has $p^{m+3}$ elements and is semilocal, but cannot be embedded in matrices over any commutative ring.
\end{thm}
Theorem 2.1 establishes that the ring  $\mathrm{End} (\mathbb{Z}_{p} \times \mathbb{Z}_{p^{m}})$ cannot be embedded in matrices over any commutative ring. Nevertheless, we can obtain a matrix representation of the elements of this ring.

Let
$$E_{p,p^m}=\{
\begin{pmatrix}a&b\\p^{m-1}c&d\end{pmatrix}|a,b,c\in\mathbb{Z}_{p},~\mathrm{and}~d\in \mathbb{Z}_{p^{m}}\},$$
where  $p^{m-1}c$ represses an element of $\mathbb{Z}_{p^{m}}$ for any $c\in\mathbb{Z}_{p}$.
\begin{lem}
The set $E_{p,p^{m}}$ is a nonccommutative unitary ring where addition is defined by
$$\begin{pmatrix}a_1&b_1\\p^{m-1}c_1&d_1\end{pmatrix}+\begin{pmatrix}a_2&b_2\\p^{m-1}c_2&d_2\end{pmatrix}=\begin{pmatrix}(a_1+a_2)~\mathrm{mod}~p&(b_1+b_2)~\mathrm{mod}~p\\p^{m-1}(c_1+c_2)~\mathrm{mod}~p^m&(d_1+d_2)~\mathrm{mod}~p^m\end{pmatrix}$$
and multiplication is defined by
$$\begin{pmatrix}a_1&b_1\\p^{m-1}c_1&d_1\end{pmatrix}\cdot\begin{pmatrix}a_2&b_2\\p^{m-1}c_2&d_2\end{pmatrix}=\begin{pmatrix}(a_1a_2)~\mathrm{mod}~p&(a_1b_2+b_1d_2)~\mathrm{mod}~p\\p^{m-1}(c_1a_2+c_2d_1)~\mathrm{mod}~p^m&(p^{m-1}c_1b_2+d_1d_2)~\mathrm{mod}~p^m\end{pmatrix}.$$
\end{lem}
\pf The proof is straightforward.
\qed

Following these results above, we now construct a ring isomorphism from $\mathrm{End} (\mathbb{Z}_{p} \times \mathbb{Z}_{p^{m}})$ onto $E_{p,p^{m}}$.
\begin{thm} With notations as above. Define the map $\varphi$ as follows
$$\varphi:~~~~~~~~\mathrm{End} (\mathbb{Z}_{p} \times \mathbb{Z}_{p^{m}})\rightarrow E_{p,p^{m}}$$
$$~~~~~~~~~~~~~~~\alpha\mapsto \varphi(\alpha)=\begin{pmatrix}a&b\\p^{m-1}c&d\end{pmatrix}$$
where $\alpha(1_p,0)=(a,p^{m-1}c)$,  $\alpha(0,1_{p^m})=(b,d)$ and $a,b,c\in\mathbb{Z}_{p}$, $d\in\mathbb{Z}_{p^{m}}$. Then  $\varphi$ is a ring isomorphism from $\mathrm{End} (\mathbb{Z}_{p} \times\mathbb{Z}_{p^{m}})$ onto $E_{p,p^{m}}$.
\end{thm}
\pf Obviously, $\varphi$ is a well-defined and an bijective mapping  from $\mathrm{End} (\mathbb{Z}_{p} \times \mathbb{Z}_{p^{m}})$ onto $E_{p,p^{m}}$.

Next, suppose that $\alpha_i\in\mathrm{End} (\mathbb{Z}_{p} \times \mathbb{Z}_{p^{m}})$ with $\alpha_{i}(1_p,0)=(a_i,p^{m-1}c_i)$, and $\alpha_i(0,1_{p^m})=(b_i,d_i)$ for $i=1,2$. Then we have  $\varphi(\alpha_i)=\begin{pmatrix}a_i&b_i\\p^{m-1}c_i&d_i\end{pmatrix}$ for $i=1,2$.
From
$$(\alpha_1+\alpha_2)(1_p,0)=(a_1,p^{m-1}c_1)+(a_2,p^{m-1}c_2)=((a_1+a_2)~\mathrm{mod}~p,p^{m-1}(c_1+c_2)~\mathrm{mod}~p^m),$$
$$(\alpha_1+\alpha_2)(0,1_{pm})=(b_1,d_1)+(b_2,d_2)=((b_1+b_2)~\mathrm{mod}~p,(d_1+d_2)~\mathrm{mod}~p^m),$$
and by the defenition of $\varphi$ and Lemma 2.2, we obtain
$$\varphi(\alpha_1+\alpha_2)=\begin{pmatrix}(a_1+a_2)~\mathrm{mod}~p&(b_1+b_2)~\mathrm{mod}~p\\p^{m-1}(c_1+c_2)~\mathrm{mod}~p^m&(d_1+d_2)~\mathrm{mod}~p^m\end{pmatrix}~~~~~~~$$
$$~~~~~~~~~=\begin{pmatrix}a_1&b_1\\p^{m-1}c_1&d_1\end{pmatrix}+\begin{pmatrix}a_2&b_2\\p^{m-1}c_2&d_2\end{pmatrix}=\varphi(\alpha_1)+\varphi(\alpha_2).$$
Moreover, by the definition of composition of endormorphism, we have
$$(\alpha_1\alpha_2)(1_p,0)=\alpha_1(a_2,p^{m-1}c_2)=\alpha_1(a_2(1_p,0)+p^{m-1}c_2(0,1_{p^m}))$$
$$=(a_2)\alpha_1(1_p,0)+p^{m-1}c_2\alpha_1(0,1_{p^m})$$$$=a_2((a_1,p^{m-1}c_1)+p^{m-1}c_2(b_1,d_1)$$
$$~~~~~~~~~~~=((a_1a_2)~\mathrm{mod}~p,p^{m-1}(a_2c_1+d_1c_2)~\mathrm{mod}~p^m)$$
and
$$(\alpha_1\alpha_2)(0,1_{p^m},)=\alpha_1(b_2,d_2)=\alpha_1(b_2(1_p,0)+d_2(0,1_{p^m}))$$
$$=b_2\alpha_1(1_p,0)+d_2\alpha_1(0,1_{p^m})~~~~~~~~~~~~~$$
$$=b_2((a_1,p^{m-1}c_1)+d_2(b_1,d_1)~~~~~~~~~~~~~~~~$$
$$~~~~~~~~~~~~~~~~~=((a_1b_2+d_2b_1)~\mathrm{mod}~p,(p^{m-1}b_2c_1+d_1d_2)~\mathrm{mod}~p^m).~~~~~~$$
From these and by Lemma 2.2, we obtain
$$\varphi(\alpha_1\alpha_2)=\begin{pmatrix}(a_1a_2)~\mathrm{mod}~p&(b_2a_1+d_1d_2)~\mathrm{mod}~p\\p^{m-1}(c_1a_2+c_2d_1)~\mathrm{mod}~p^m&(p^{m-1}c_1b_2+d_1d_2)~\mathrm{mod}~p^m\end{pmatrix}~~~~~~~$$
$$~~~~~~~~~=\begin{pmatrix}a_1&b_1\\p^{m-1}c_1&d_1\end{pmatrix}\cdot \begin{pmatrix}a_2&b_2\\p^{m-1}c_2&d_2\end{pmatrix}=\varphi(\alpha_1)\cdot\varphi(\alpha_2).~~~~~~~~~~~~~~~~~~~~~$$
This means that $\varphi$ is a ring isomorphism from $\mathrm{End} (\mathbb{Z}_{p} \times\mathbb{Z}_{p^{m}})$ onto $E_{p,p^{m}}$.
\qed
\begin{rem}$| \mathrm{End} (\mathbb{Z}_{p} \times\mathbb{Z}_{p^{m}})|=| E_{p,p^{m}}|.$
\end{rem}

From  now on, we identity the elements of  $\mathrm{End} (\mathbb{Z}_{p} \times\mathbb{Z}_{p^{m}})$  with the elements of  $E_{p,p^{m}} $, and the arithmetic of $\mathrm{End} (\mathbb{Z}_{p} \times\mathbb{Z}_{p^{m}})$ with the  arithmetic of  $E_{p,p^{m}}$.

 \section{The arithmetic of invertible elements in the ring $E_{p,p^{m}}$}
In this section we present our contributions, i.e., we first provide the way to compute $-d$ and $d^{-1}$ using arithmetic  in $\mathbb{Z}_{p}$ and $\mathbb{Z}_{p^{m}}$. We then characterize invertible elements in the ring $E_{p,p^{m}}$.

As usual, if $a,b\in\mathbb{Z}$, with $b\neq0$, we denote by $\lfloor\frac{a}{b}\rfloor$ and $a~\mathrm{mod}~b$ the quotient and the remainder of the division of $a$ by $b$, respectively.

We first observe four helpful lemmas.
\begin{lem}
 Suppose that $d_i=p^{m-1}u_{m-1}^{(i)}+\cdots+pu_{1}^{(i)}+u_{0}^{(i)}$with  $u_{m-1}^{(i)},\ldots,u_{1}^{(i)},u_{0}^{(i)}\in \mathbb{Z}_{p}$, for $i=1,2$. Let
$$
v_0=(u_{0}^{(1)}+u_{0}^{(2)}) ~\mathrm{mod}~p^{m},~~~~~~~~~~~~~~~~~~~~$$
$$v_{1}=(u_{1}^{(1)}+u_{1}^{(2)}+\lfloor\frac{v_{0}}{p}\rfloor)~\mathrm{mod}~p^{m},~~~~~~~~~~~$$
$$v_2=(u_{2}^{(1)}+u_{2}^{(2)}+\lfloor\frac{v_{1}}{p}\rfloor)~\mathrm{mod}~p^{m},~~~~~~~~~~~$$
$$\cdots\\$$
$$v_{m-2}=(u_{m-2}^{(1)}+u_{m-2}^{(2)}+\lfloor\frac{v_{m-3}}{p}\rfloor)~\mathrm{mod}~p^{m},$$
$$v_{m-1}=(u_{m-1}^{(1)}+u_{m-1}^{(2)}+\lfloor\frac{v_{m-2}}{p}\rfloor)~\mathrm{mod}~p^{m}.$$
Then $d_1+d_2=p^{m-1}(v_{m-1}~\mathrm{mod}~p)+\cdots+p(v_{1}~\mathrm{mod}~p)+(v_{0}~\mathrm{mod}~p).$
\end{lem}
\pf Obviously,
$$d_1+d_2=[p^{m-1}(u_{m-1}^{(1)}+u_{m-1}^{(2)})+\cdots+p(u_{1}^{(1)}+u_{1}^{(2)})+(u_{0}^{(1)}+u_{0}^{(2)})]~\mathrm{mod}~p^m.$$
According to the definition of $v_i$ we have
$$d_1+d_2=p^{m-1}[(u_{m-1}^{(1)}+u_{m-1}^{(2)}+\lfloor\frac{v_{m-2}}{p}\rfloor)~\mathrm{mod}~p]+$$
$$~~~~~~~~~~~~~~~~~~~~~~~~~~~~~~\cdots+p[(u_{1}^{(1)}+u_{1}^{(2)}+\lfloor\frac{v_{0}}{p}\rfloor)~\mathrm{mod}~p]+(u_{0}^{(1)}+u_{0}^{(2)})~\mathrm{mod}~p$$ $$~~~~~~~~~~~~~~~~~~~~~~~~~~=p^{m-1}(v_{m-1}~\mathrm{mod}~p)+\cdots+p(v_{1}~\mathrm{mod}~p)+(v_{0}~\mathrm{mod}~p).$$
\qed
\begin{lem}
 Suppose that $d=p^{m-1}u_{m-1}+\cdots+pu_{1}+u_{0}$  with  $u_{m-1},\ldots,u_{1},u_{0}\in \mathbb{Z}_{p}$ . Let
$$v_0=-u_{0}~\mathrm{mod}~p,~~~~~~~~~~~~~~~~~$$
$$w_{0}=(u_{0}-u_{0})~\mathrm{mod}~p^{m}~~~~~~~~~~~~~~$$
$$v_{1}=(-u_{1}-\lfloor\frac{w_{0}}{p}\rfloor )~\mathrm{mod}~p,~~~~~~~$$
$$w_{1}=(u_{1} +v_1+\lfloor\frac{w_{0}}{p}\rfloor) \mathrm{mod}~p^{m}~~~~~~~~~~$$
$$v_{2}=(-u_{2}-\lfloor\frac{w_{1}}{p}\rfloor )~\mathrm{mod}~p,~~~~~~~$$
$$w_{2}=(u_{2}+v_{2})+\lfloor\frac{w_{1}}{p}\rfloor ~\mathrm{mod}~p^{m}~~~~~~~~~~$$
$$\cdots\\$$
$$v_{m-2}=(-u_{m-2}-\lfloor\frac{w_{m-3}}{p}\rfloor )~\mathrm{mod}~p,~~~~~~~$$
$$w_{m-2}=(u_{m-2}+v_{m-2}+\lfloor\frac{w_{m-3}}{p}\rfloor, ~\mathrm{mod}~p^{m}~$$
$$v_{m-1}=(-u_{m-1}-\lfloor\frac{w_{m-2}}{p}\rfloor )~\mathrm{mod}~p.~~~~~~~$$
Then  $-d=p^{m-1}v_{m-1}+\cdots+pv_{1}+v_{0}.$
\end{lem}
\pf In light of Lemma 3.1,
$$d-d=p^{m-1}(s_{m-1}~\mathrm{mod}~p)+\cdots+p(s_{1}~\mathrm{mod}~p)+(s_{0}~\mathrm{mod}~p),$$
where
$$s_0=(u_0-u_0)~\mathrm{mod}~p^{m},$$
$$s_1=(u_1+[(-u_1-\lfloor\frac{w_0}{p}\rfloor)~\mathrm{mod}~p]+\lfloor\frac{w_0}{p}\rfloor)~\mathrm{mod}~p^{m},$$
$$s_2=(u_2+[(-u_2-\lfloor\frac{w_1}{p}\rfloor)~\mathrm{mod}~p]+\lfloor\frac{w_1}{p}\rfloor)~\mathrm{mod}~p^{m},$$
$$\cdots$$
 $$s_{m-1}=(u_{m-1}+[(-u_{m-1}-\lfloor\frac{w_{m-2}}{p}\rfloor)~\mathrm{mod}~p]+\lfloor\frac{w_{m-2}}{p}\rfloor)~\mathrm{mod}~p^{m}.$$
we have
$$s_0=(u_0-u_0)~\mathrm{mod}~p=0,$$
$$s_1=[u_1+((-u_1-\lfloor\frac{w_0}{p}\rfloor)~\mathrm{mod}~p)+\lfloor\frac{w_0}{p}\rfloor]~\mathrm{mod}~p=0,$$
$$s_2=[u_2+((-u_2-\lfloor\frac{w_1}{p}\rfloor)~\mathrm{mod}~p)+\lfloor\frac{w_1}{p}\rfloor]~\mathrm{mod}~p=0,$$
 $$\cdots$$
$$s_{m-1}=[u_{m-1}+((-u_{m-1}-\lfloor\frac{w_{m-2}}{p}\rfloor)~\mathrm{mod}~p)+\lfloor\frac{w_{m-2}}{p}\rfloor] ~\mathrm{mod}~p=0.$$
Therefore,
$$d-d=p^{m-1}(s_{m-1}~\mathrm{mod}~p)+\cdots+p(s_{1}~\mathrm{mod}~p)+(s_{0}~\mathrm{mod}~p)=0.$$
This proves the expected result.
\qed
\begin{lem}
Suppose that $d_i=p^{m-1}u_{m-1}^{(i)}+\cdots+pu_{1}^{(i)}+u_{0}^{(i)}$ with $u_{m-1}^{(i)},\ldots,u_{1}^{(i)},u_{0}^{(i)}\in \mathbb{Z}_{p},$  for $i=1,2$. Let
$$w_0=u_{0}^{(1)}u_{0}^{(2)}~\mathrm{mod}~p^{m},~~~~~~~~~~~~~~~~~~~~~~~~~~$$
$$w_{1}=(\sum_{i+j=1}u_{i}^{(1)}u_{j}^{(2)}+\lfloor\frac{w_{0}}{p}\rfloor)~\mathrm{mod}~p^{m},~~~~~~~~~~$$
$$w_2=(\sum_{i+j=2}u_{i}^{(1)}u_{j}^{(2)}+\lfloor\frac{w_{1}}{p}\rfloor)~\mathrm{mod}~p^{m},~~~~~~~~~~$$
$$\cdots\\$$
$$w_{m-2}=(\sum_{i+j=m-2}u_{i}^{(1)}u_{j}^{(2)}+\lfloor\frac{w_{m-3}}{p}\rfloor)~\mathrm{mod}~p^{m},$$
$$w_{m-1}=(\sum_{i+j=m-1}u_{i}^{(1)}u_{j}^{(2)}+\lfloor\frac{w_{m-1}}{p}\rfloor)~\mathrm{mod}~p^{m}.$$
Then $d_1d_2=p^{m-1}(w_{m-1}~\mathrm{mod}~p)+\cdots+p(w_{1}~\mathrm{mod}~p)+(w_{0}~\mathrm{mod}~p).$
\end{lem}
\pf The proof is similar to Lemma 3.1.
\qed
\begin{lem}
 Suppose that $d=p^{m-1}u_{m-1}+\cdots+pu_{1}+u_{0}$ with  $u_{m-1},\ldots,u_{1},u_{0}\in \mathbb{Z}_{p}$ . If $u_{0}\neq 0$, then  $d$ is invertible in $\mathrm{Z_{p^{m}}}$. Further,  write
$$s_0=u_{0}^{-1}~\mathrm{mod}~p~~~~~~~~~~~~~~~~~,$$
$$w_{0}=u_{0}u_{0}^{-1}~\mathrm{mod}~p^{m}~~~~~~~~~~~~~~$$
$$s_{1}=-u_{0}^{-1}(u_{1},1)(s_0,\lfloor\frac{w_{0}}{p}\rfloor )^{T}~\mathrm{mod}~p,~~~~~~~$$
$$w_{1}=(u_{1},u_{0})(s_0,s_1)^{T}+\lfloor\frac{w_{0}}{p}\rfloor~ \mathrm{mod}~p^{m}~~~~~~~~~~$$
$$s_{2}=-u_{0}^{-1}(u_2,u_{1},1)(s_0,s_1,\lfloor\frac{w_{1}}{p}\rfloor )^{T}~\mathrm{mod}~p,~~~~~~~$$
$$w_{2}=(u_2,u_{1},u_{0})(s_0,s_1,s_2)^{T}+\lfloor\frac{w_{1}}{p}\rfloor ~\mathrm{mod}~p^{m}~~~~~~~~~~$$
in general, for $3 \leq k \leq m-1$
$$s_{k}=-u_{0}^{-1}((u_k,\ldots,u_{1},1)(s_0,\ldots,s_{k-1},\lfloor\frac{w_{k-1}}{p}\rfloor )^{T}~\mathrm{mod}~p,$$
and,
$$w_{k}=(u_k,\ldots,u_{0})(s_0,\ldots,s_{k-1})^{T}+\lfloor\frac{w_{k-1}}{p}\rfloor ~\mathrm{mod}~p^{m}.$$
where $\mathrm{T}$ denotes transpose of matrix.
Then  $d^{-1}=p^{m-1}s_{m-1}+\cdots+ps_{1}+s_{0}.$
\end{lem}

\pf Since $0\neq u_{0}\in\mathbb{Z}_p$, there is a unique inverse $u_{0}^{-1}\in\mathbb{Z}_p$ of $u_{0}$ such that $u_{0}u_{0}^{-1}=1$ in $\mathbb{Z}_{p}$.
Take, $a=p^{m-1}t_{m-1}+\cdots+pt_{1}+t_{0}$ with $t_{m-1},\cdots,t_{1},t_{0}\in\mathrm{Z}_{p^{m}}$, then by Lemma 3.2, we have
$ad=p^{m-1}(h_{m-1}~\mathrm{mod}~p)+\cdots+p(h_{1}~\mathrm{mod}~p)+(h_{0}~\mathrm{mod}~p),$
where
$$
h_0=u_{0}t_{0},~~~~~~~~~~~~~~~~~~~~~~~~~~$$
$$h_{1}=\sum_{i+j=1}u_{i}t_{j}+\lfloor\frac{h_{0}}{p}\rfloor,~~~~~~~~~~$$
$$h_2=\sum_{i+j=2}u_{i}t_{j}+\lfloor\frac{h_{1}}{p}\rfloor,~~~~~~~~~~$$
$$\cdots\\$$
$$h_{m-2}=\sum_{i+j=m-2}u_{i}t_{j}+\lfloor\frac{h_{m-3}}{p}\rfloor,$$
$$h_{m-1}=\sum_{i+j=m-1}u_{i}t_{j}+\lfloor\frac{h_{m-2}}{p}\rfloor.$$

In order to get $ad=1$ in $\mathrm{Z_p^{m}}$, first taking $t_0=s_{0}=u_{0}^{-1}$, we have $h_0=w_0=u_{0}u_{0}^{-1}$=1 in $\mathbb{Z}_{p}$.

Next, taking, $t_1=s_{1}$, we have $h_1=w_1=0$ in $\mathbb{Z}_p$,  but $h_1=w_1$ in $\mathbb{Z}_p^{m}$.

Continuing this process, taking $t_k=s_k$, we have $h_k=w_k=0$ in $\mathbb{Z}_p$, but $h_k=w_k$ in $\mathbb{Z}_p^{m}$, for $k=2,3,\ldots,m-1$. This means that
$$a=p^{m-1}s_{m-1}+\cdots+ps_{1}+s_{0}.$$
Therefore $ad=1$ in $\mathbb{Z}_p^{m}$, i.e., $d^{-1}=p^{m-1}s_{m-1}+\cdots+ps_{1}+s_{0}.$
\qed

In the following we establish a characterization of the invertible element of $E_{p,p^m}$. First, by  Lemma 3.1 and Lemma 3.2, we obtain the following characterization of addition and multiplication in $Z_{p,p^{m}}$ in terms of the arithmetic of $\mathbb{Z_p}$ and $\mathbb{Z}_{p^m}$.
\begin{cor}
Assume that all hypothesis of Lemma 3.1 and Lemma 3.3 are valid. Let
\[\begin{matrix}
A_1=\begin{pmatrix}a_1&b_1\\p^{m-1}c_1&p^{m-1}u_{m-1}^{(1)}+\ldots+pu_1^{(1)}+u_0^{(1)}\end{pmatrix}
\end{matrix}\]
and
\[\begin{matrix}
A_2=\begin{pmatrix}a_2&b_2\\p^{m-1}c_2&p^{m-1}u_{m-1}^{(2)}+\ldots+pu_1^{(2)}+u_0^{(2)}\end{pmatrix}
\end{matrix}\]
be two elements of $E_{p,p^m}$. Then
\[\begin{matrix}
A_1+A_2=\begin{pmatrix}(a_1+a_2)~\mathrm{mod}~p&(b_1+b_2))~\mathrm{mod}~p\\p^{m-1}((c_1+c_2)~\mathrm{mod}~p)&p^{m-1}(v_{m-1}~\mathrm{mod}~p)+\ldots+p(v_1~\mathrm{mod}~p)+v_o~\mathrm{mod}~p\end{pmatrix},
\end{matrix}\]
and
\[\begin{matrix}
A_1A_2=\begin{pmatrix}(a_1a_2)~\mathrm{mod}~p&(a_1b_2+b_1u_0^{(2)}))~\mathrm{mod}~p\\p^{m-1}((c_1a_2+u_0^{(1)}c_2)~\mathrm{mod}~p)&\delta\end{pmatrix}
\end{matrix}\]
where $\delta= p^{m-1}((c_1b_2+w_{m-1})~\mathrm{mod}~p)+\ldots+pw_1+w_0$.
\end{cor}
Next, combining Lemma 3.4 and Corollary 3.5, we can now to  establish a characterization of the invertible element of $E_{p,p^{m}}$.
\begin{thm}
Suppose that
\[\begin{matrix}
M=\begin{pmatrix}a&b\\p^{m-1}c&p^{m-1}u_{m-1}+\ldots+pu_1+u_0\end{pmatrix}\in E_{p,p^{m}},
\end{matrix}\]
where $a,b,c,u_0,u_1,\ldots,u_{m-1}\in \mathbb{Z}_p$. Then $M$ is invertible if and only if $a\neq 0$ and $ u_0\neq 0$, and in this case the inverse $M^{-1}$ of $M$ in $E_{p,p^{m}}$ is given by
$$\begin{matrix}
M^{-1}=\begin{pmatrix}a^{-1}\mathrm{mod}p&(-a^{-1}bu_0^{-1})\mathrm{mod}p\\p^{m-1}((-a^{-1}cu_0^{-1})\mathrm{mod}p)&p^{m-1}[(ca^{-1}bu_0^{-2}+s_{m-1})\mathrm{mod}p]+p^{m-2}s_{m-2}+\ldots+ps_1+s_0\end{pmatrix},
\end{matrix}$$
where $p^{m-1}s_{m-1}+\ldots+ps_1+s_0$ is invertible of $p^{m-1}u_{m-1}+\ldots+pu_1+u_0$.
\end{thm}
\pf If $M$ is invertible, the there exists
\[\begin{matrix}
N=\begin{pmatrix}x&y\\p^{m-1}z&p^{m-1}t_{m-1}+\ldots+pt_1+t_0\end{pmatrix}\in E_{p,p^{m}},
\end{matrix}\]
with $x,y,z,t_0,t_1,\ldots,t_{m-1}\in \mathbb{Z}_p$, such that $MN=I$, where
$\begin{matrix}
I=\begin{pmatrix}1&0\\0&1\end{pmatrix}
\end{matrix}$.  In light of Corollary 3.4, we obtain
$$ 1=(ax)~\mathrm{mod}~p~\mathrm{and}~1=(u_0t_0)~\mathrm{mod}~p.$$
Therefore $a\neq 0$ and $ u_0\neq 0$.

Reciprocally, let $a\neq 0$ and $ u_0\neq 0$. Then there exist $a^{-1},u_0^{-1}\in\mathbb{Z}_p$. Taking
\[\begin{matrix}
N=\begin{pmatrix}a^{-1}\mathrm{mod}p&(-a^{-1}bu_0^{-1})\mathrm{mod}p\\p^{m-1}((-a^{-1}cu_0^{-1})\mathrm{mod}p)&p^{m-1}[(ca^{-1}bu_0^{-2}+s_{m-1})\mathrm{mod}p]+p^{m-2}s_{m-2}+\ldots+ps_1+s_0\end{pmatrix},
\end{matrix}\]
where $p^{m-1}s_{m-1}+\ldots+ps_1+s_0$ is invertible of $p^{m-1}u_{m-1}+\ldots+pu_1+u_0$. From Corollary 3.5, we have that $\begin{matrix}
I=\begin{pmatrix}h&r\\p^{m-1}l&t\end{pmatrix},
\end{matrix}$
where
$$h=(aa^{-1}~\mathrm{mod}~p=1,$$
$$r=(a(-a^{-1}bu_0^{-1})+bu_0^{-1})~\mathrm{mod}~p=0,$$
$$l=(ca^{-1}-u_0^{-1}ca^{-1}u_0)~\mathrm{mod}~p=0,$$
and
$$t=p^{m-1}c(-a^{-1}bu_0^{-1})+[p^{m-1}(ca^{-1}bu_0^{-2}+s_{m-1})+$$
$$~~~~~~~~~~~~~~~~~~~~p^{m-2}s_{m-1}+\ldots+ps_1+s_0][p^{m-1}u_{m-1}+\ldots+pu_1+u_0]=1.$$
And consequently, $MN=I$.

Similar to prove $NM=I$. therefore, $M$ is invertible in $E_{p,p^{m}}$ and $M^{-1}=N$.
\qed

The following result gives the number of invertible elements of $E_{p,p^{m}}$.
\begin{cor}
Let $U_{p,p^{m}}$ be the set of all invertible elements of $E_{p,p^{m}}$. Then $|U_{p,p^{m}}|=p^{m+1}(p-1)^{2}$.
\end{cor}
\pf By Theorem 3.6, the number of invertible elements of $E_{p,p^{m}}$ equals $|U_{p,p^{m}}|=|\mathbb{Z}_{p}|^{m-1}|\mathbb{Z}_{p}-0|^{2}=p^{m+1}(p-1)^{2}.$
\qed

Note that $\frac{U_{p,p^{m}}}{E_{p,p^{m}}}=\frac{p^{m+1}(p-1)^{2}}{p^{m+3}}=(1-\frac{1}{p})^{2}\approx1$. This means that for large values of $p$, almost all the elements of $E_{p,p^{m}}$ are invertible.

\section{The minimal polynomial of each element in $E_{p,p^{m}}$}
In this section, we introduce a concept of the  minimal polynomial for each element in $E_{p,p^{m}}$.
We start by showing Theorem 4.1.
\begin{thm} Suppose that $\begin{matrix}
A=\begin{pmatrix}a&b\\p^{m-1}c&d\end{pmatrix}
\end{matrix}\in E_{p,p^{m}}$ with $a,b,c\in\mathbb{Z}_{p}$ and $d=p^{m-1}u_{m-1}+\ldots+pu_1+u_0$ where $u_0,u_1,\ldots,u_{m-1}\in\mathbb{Z}_{p}$.  Let $$r=-(a+d)~\mathrm{mod}~p^{m},~~~~~~~s=(ad-p^{m-1}bc)~\mathrm{mod}~p^{m},$$
and denote $f_A(x)=x^{2}+rx+s\in \mathbb{Z}[x]$.  Then $f_A(A)=A^{2}+rA+sI=0.$
\end{thm}
\pf By assumption, we have
\[\begin{matrix}
M^2=\begin{pmatrix}a^2~\mathrm{mod}~p&(ab+bu_0)~\mathrm{mod}~p\\p^{m-1}[(ca+cu_0)~\mathrm{mod}~p]&(p^{m-1}cb+d^2))~\mathrm{mod}~p^{m}\end{pmatrix}\in E_{p,p^{m}},
\end{matrix}\]
\[\begin{matrix}
rM=\begin{pmatrix}[(-a^2-ad)~\mathrm{mod}~p^{m}]~\mathrm{mod}~p&[(-ab+bd)~\mathrm{mod}~p^{m}]~\mathrm{mod}~p\\-p^{m-1}[(ca+cu_0)~\mathrm{mod}~p^{m}]~\mathrm{mod}~p]&-((a+d)~\mathrm{mod}~p^{m})\cdot(d~\mathrm{mod}~p^{m})\end{pmatrix}\\~~~~~~~~~~~~~~=\begin{pmatrix}(-a^2-au_0)~\mathrm{mod}~p&(-ab+bu_0)~\mathrm{mod}~p\\-p^{m-1}[(ca+cu_0)~\mathrm{mod}~p]&-(ad)~\mathrm{mod}~p^{m}-d^2~\mathrm{mod}~p^{m}\end{pmatrix},~~~~~~~~~~~~~~~~~~~
\end{matrix}\]
and
\[\begin{matrix}
sI=\begin{pmatrix}[(ad-p^{m-1}bc)~\mathrm{mod}~p^{m}]~\mathrm{mod}~p&0\\0&(ad-p^{m-1}bc)~\mathrm{mod}~p^{m}\end{pmatrix}.
\end{matrix}\]
Therefore $f_A(A)=A^{2}+rA+sI=0.$
\qed
\begin{defi}
For $A\in E_{p,p^{m}}$, we call monic polynomial $f_{A}(x)\in \mathbb{Z}[x]$  the minimal polynomial of $A$ if $f_{A}(A)=0$ in $E_{p,p^{m}}$ and $g(A)=0$ if and only if $f_{A}(x)$ divides $g(x)$ in $\mathbb{Z}[x]$.
\end{defi}
\begin{cor} In the notation of Theorem 4.1, $f_A(x)=x^{2}+rx+s$ is the minimal polynomial of $A$.
\end{cor}

In $\mathbb{Z}[x]$ or  $\mathbb{Z}_{p^m}[x]$, for any $g(x)\in \mathbb{Z}[x]$, by the polynomial division with remainder algorithm over $\mathbb{Z}[x]$, we obtain
$$g(x)=q(x)f_A(x)+r(x),$$
where $q(x),r(x)\in \mathbb{Z}[x]$, and $r(x)=0$ or $deg(r(x))<deg(f_A(x))$. Hence
$$g(A)=q(A)f_A(A)+r(A)=r(A).$$
This means that we can simplify the computations of $g(A)$ for any $g(x)\in \mathbb{Z}[x]$ and $M^{-1}$ in $E_{p,p^{m}}$.

\begin{cor} With notations of Theorem 4.1.  If $a\neq0$ and $u_0\neq0$, then $A^{-1}=-s^{-1}(A+rI).$
\end{cor}
\pf Since $A^{2}+rA+sI=0$,  we have $(A+rI)A=(-s)I$. By $a\neq0$ and $u_0\neq0$. Therefore,
$$A^{-1}=(-s)^{-1}(A+rI)=-s^{-1}(A+rI).$$
\qed

We finish this section by showing how Corollary 4.4 works for computing invertible elements in $E_{p,p^{m}}$.

\begin{exam} Suppose that $\begin{matrix}
A=\begin{pmatrix}2&3\\5^{2}\cdot3&5^2\cdot2+5\cdot3+2\end{pmatrix}
\end{matrix}\in E_{5,5^{3}}$. Then by Theorem 4.1 we have
$$r=-(2+5^2\cdot2+5\cdot3+2)~\mathrm{mod}~5^3=56~\mathrm{mod}~5^3=5^2\cdot2+5\cdot1+1,$$
and
$$s=[2(5^2\cdot2+5\cdot3+2)-5^2\cdot3\cdot3]~\mathrm{mod}~5^3=34~\mathrm{mod}~5^3=5^2\cdot1+5\cdot1+4.$$
In light of  Lemma 3.3, we obtain
$$s^{-1}=5^2\cdot4+\cdot2+4,$$
hence,$-s^{-1}=5\cdot2+1$.

Now, by Corollary 4.4, we have
$$\begin{matrix}A^{-1}=(5\cdot2+1)[\begin{pmatrix}2&3\\5^{2}\cdot3&5^2\cdot2+5\cdot3+2\end{pmatrix}+\begin{pmatrix}4&0\\0&5^2\cdot2+5\cdot1+1\end{pmatrix}]\\~~~~~~~~~~~~~~=(5\cdot2+1)\begin{pmatrix}3&3\\5^2\cdot3&5^2\cdot4+5\cdot4+3\end{pmatrix}=\begin{pmatrix}3&3\\5^2\cdot3&5^2\cdot4+3\end{pmatrix}.~~~~~~~~~~~~
\end{matrix}$$
\end{exam}

\end{document}